# Rational Orthonormal Littlewood-Paley Wavelet Basis of $L^2(R)$

Sun Changping, *Member*, *IEEE*

Department of Artificial Intelligence, School of Electrical Engineering and Automation, Tianjin University of Technology, Tianjin, 300384, China. E-mail: sunchangping@tjut.edu.cn

*Abstract*—In this letter, first, we prove that the orthonormal basis of rational Littlewood-Paley wavelet with rational dilation factor M=p/q first proposed by Auscher does not hold for all rational numbers. It does not hold if q is not equal to 1. In other words, it is not an orthonormal basis if the rational dilation factor M is not an integer. Then, to make up for the shortcoming of the rational Littlewood-Paley wavelet proposed by Auscher, a new orthonormal basis of rational Littlewood-Paley wavelet with rational dilation factor M=p/q is proposed, which holds for all rational numbers. Finally, by means of sampling theorem for bandpass signals, it is proved completely that the new rational Littlewood-Paley wavelet family $\left\{ \psi_{j,n}^m(x) = M^{-\frac{j}{2}} \psi^m(M^{-j}x - nq) \right\}_{j,n \in \mathbb{Z},\ 1 \leq m \leq p-q}$ is an orthonormal wavelet basis of $L^2(R)$.

## I. Introduction

In [1][2], the author first proposed the concept of rational multiresolution analysis. It could provide a better separation of signal components than the dyadic one [3]. In the rational multiresolution analysis, the construction of rational wavelets is vital. The orthonormal basis of rational Littlewood-Paley wavelet was first proposed in [1]. Although the form of it in the frequency domain is relatively simple, it plays an important role in rational multiresolution analysis. For example, it can be used to verify the validity of the algorithm of rational wavelet decomposition and reconstruction. If the basis of rational Littlewood-Paley wavelet is not orthonormal, this can lead to wrong results. For example, if it used to verify the validity of the algorithm of rational wavelet decomposition and reconstruction, the original signal can not be exactly synthesized based on the rational wavelet decomposition coefficients, although the algorithm of rational wavelet decomposition and synthesis is completely correct in theory. A pyramidal algorithm for fast rational orthogonal wavelet transform was proposed in [3]. The condition characterizing completeness of orthonormal wavelet systems with arbitrary real dilation factor was given in [4]. In [5], the authors gave a complete characterization of tight frames of orthonormal wavelets with arbitrary real dilation factor. In [6], a perfect reconstruction condition was given for the orthonormal wavelet with rational dilation factor. An overcomplete discrete wavelet transform with rational dilation factor was developed in [7]. In [8], a discrete wavelet transform fast algorithm with dilation factor 3/2 is proposed, which overcomes the drawback of frequency distortion in the high frequency subband in the process of decomposition of Mallat discrete wavelet transform (DWT). In addition, wavelet transform with rational dilation factor has been applied to fault diagnosis [9], automatic imagined speech recognition [10], and automatic detection of heart valve disorders [11]. Although rational wavelets are widely researched, to the best of the author's knowledge, the orthonormal rational Littlewood-Paley wavelet cited in the current literature is still the form proposed by Auscher, and there are no proofs of the rational Littlewood-Paley wavelet family being an orthonormal wavelet basis of $L^2(R)$. Auscher originally intended to give the orthonormal basis of rational Littlewood-Paley wavelet which holds for all rational dilation factor, but it is a pity that Auscher's orthonormal rational Littlewood-Paley wavelet with rational dilation factor M=p/q does not actually hold for all rational numbers. It does not hold if q is not equal to 1. In other words, it is not an orthonormal basis if the rational dilation factor M is not an integer.

The rest of this study is organized as follows. The limitation of Auscher's rational Littlewood-Paley wavelet is pointed out, and the proof is given in section II. In section III, a new orthonormal basis of rational Littlewood-Paley wavelet with rational dilation factor M=p/q is proposed. In section IV, it is proved completely that the new rational Littlewood-Paley wavelet family $\left\{\psi_{j,n}^m(x) = M^{-\frac{j}{2}} \psi^m(M^{-j}x - nq)\right\}_{j,n \in \mathbb{Z},\ 1 \leq m \leq p-q}$ is an orthonormal wavelet basis of $L^2(R)$. Finally, the conclusion is given in section V.

## II. The Limitation of Auscher's Rational Littlewood-Paley Wavelet

The orthonormal basis of rational Littlewood-Paley wavelet with rational dilation factor M=p/q was first proposed by Auscher in [1]. However, Auscher originally intended to give the orthonormal basis of rational Littlewood-Paley wavelet which holds for all rational numbers, but it is a pity that Auscher's rational Littlewood-Paley wavelet with rational dilation factor M=p/q does not actually hold for all rational numbers, that is, it does not hold if q is not equal to 1.

**Theorem 1:** The basis of rational Littlewood-Paley wavelet with rational dilation factor M=p/q proposed by Auscher is not orthonormal if q is not equal to 1.

**Proof:** In the frequency domain, the rational Littlewood-Paley wavelet proposed by Auscher is given in formula (1) [1].

$$\hat{\psi}^m(M\omega) = \begin{cases} 1 & \text{if } (q+m-1)\frac{\pi}{p} \leq |\omega| < (q+m)\frac{\pi}{p} \\ 0 & \text{elsewhere} \end{cases} \tag{1}$$

Where $M = \frac{p}{q}$, $1 \leq m \leq p-q$, $p$ and $q$ relatively prime, $p > q \geq 1$.

By computing

$$\|\psi_{j,n}^m(x)\|^2$$

$$= \langle \psi_{j,n}^m(x),\ \psi_{j,n}^m(x)\rangle$$

$$= \frac{1}{2\pi}\left\langle \widehat{\psi}_{j,n}^m(\omega),\ \widehat{\psi}_{j,n}^m(\omega)\right\rangle$$

$$= \frac{1}{2\pi}\int_{-\infty}^{+\infty} M^{\frac{j}{2}}\widehat{\psi}^m(M^j\omega)e^{-inM^jq\omega}\ \overline{M^{\frac{j}{2}}\widehat{\psi}^m(M^j\omega)e^{-inM^jq\omega}}\ d\omega$$

$$= \frac{M^j}{2\pi}\left[\int_{M^{-j}(m+q-1)\frac{\pi}{q}}^{M^{-j}(m+q)\frac{\pi}{q}} 1\ d\omega + \int_{-M^{-j}(m+q)\frac{\pi}{q}}^{-M^{-j}(m+q-1)\frac{\pi}{q}} 1\ d\omega\right]$$

$$= \frac{M^j}{2\pi}\left[M^{-j}\frac{\pi}{q} + M^{-j}\frac{\pi}{q}\right]$$

$$= \frac{1}{q} \tag{2}$$

It can be found that if q is not equal to 1, the norm of $\psi_{j,n}^m(x)$ is not equal to 1. So the basis of rational Littlewood-Paley wavelet with rational dilation factor M=p/q proposed by Auscher is not orthonormal if q is not equal to 1. □

### III. The New Rational Littlewood-Paley Wavelet Basis

In this section, the new rational Littlewood-Paley wavelet is given in frequency domain and time domain, respectively.

In the frequency domain, the new rational Littlewood-Paley wavelet is given as follows.

$$\widehat{\psi}^m(M\omega) = \begin{cases} \sqrt{q} & \text{if } (q+m-1)\frac{\pi}{p} \leq |\omega| < (q+m)\frac{\pi}{p} \\ 0 & \text{elsewhere} \end{cases} \tag{3}$$

Where $M = \frac{p}{q}$, $p$ and $q$ relatively prime, $p > q \geq 1$.

Based on the formula (3), the formulas of $\widehat{\psi}^m(\omega)$ and $\widehat{\psi}^m(M^j\omega)$ are given as follows.

$$\widehat{\psi}^m(\omega) = \begin{cases} \sqrt{q} & \text{if } (q+m-1)\frac{\pi}{q} \leq |\omega| < (q+m)\frac{\pi}{q} \\ 0 & \text{elsewhere} \end{cases} \tag{4}$$

$$\widehat{\psi}^m(M^j\omega) = \begin{cases} \sqrt{q} & \text{if } M^{-j}(q+m-1)\frac{\pi}{q} \leq |\omega| < M^{-j}(q+m)\frac{\pi}{q} \\ 0 & \text{elsewhere} \end{cases} \tag{5}$$

According to the formula (4), the new rational Littlewood-Paley wavelet in time domain is given as follows:

$$\psi^m(t) = \sqrt{q}\left[\frac{\sin(q+m)\frac{\pi}{q}t}{\pi t} - \frac{\sin(q+m-1)\frac{\pi}{q}t}{\pi t}\right] \tag{6}$$

Based on the formula (3), the new rational Littlewood-Paley wavelets $\hat{\psi}^m(M\omega)$ for $M=\frac{5}{3}$ are shown in Fig.1 and Fig.2, respectively, here *m*=1, 2.

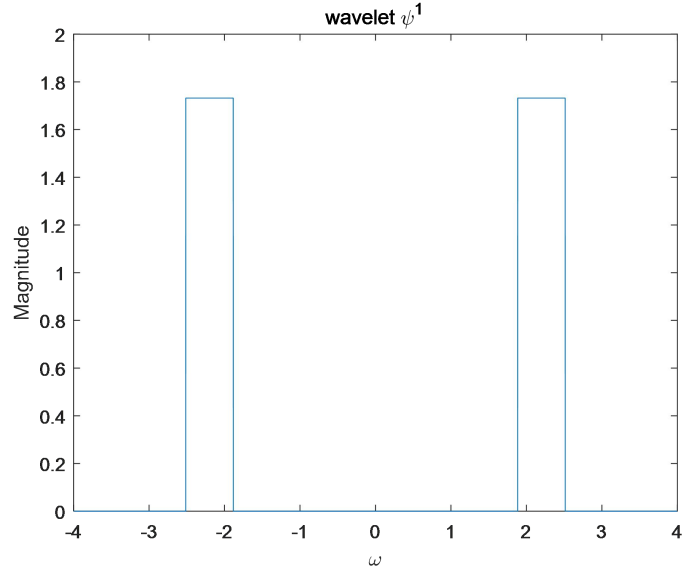

Fig. 1. The new rational Littlewood-Paley wavelet $\hat{\psi}^m(M\omega)$ for $M=\frac{5}{3}$, $m=1$

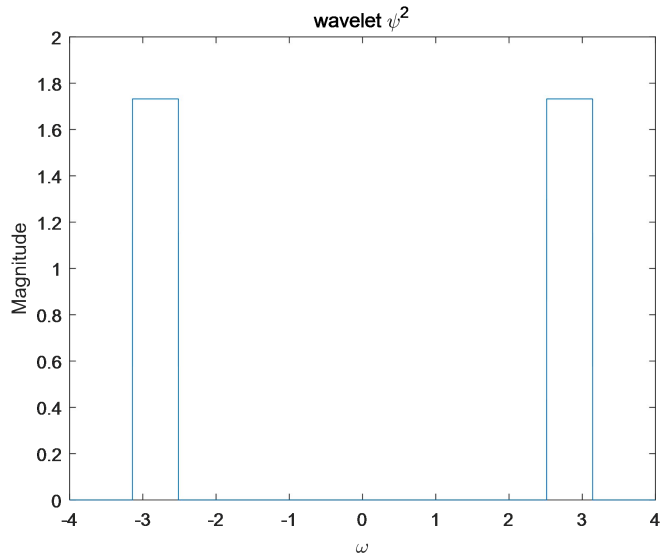

Fig. 2. The new rational Littlewood-Paley wavelet $\hat{\psi}^m(M\omega)$ for $M=\frac{5}{3}$, $m=2$

# IV. The Proof of the New Rational Littlewood-Paley Wavelet Family is an Orthonormal Wavelet Basis for $L^2(R)$

In this section, the complete proof of the new rational Littlewood-Paley wavelet set $\left\{\psi_{j,n}^m(x) = M^{-\frac{j}{2}}\psi^m(M^{-j}x - nq)\right\}_{j,n \in \mathbb{Z},\ 1 \leq m \leq p-q}$ is an orthonormal wavelet basis for $L^2(R)$ is given.

**Theorem 2:** In the frequency domain, if the rational Littlewood-Paley wavelet basis $\hat{\psi}^m(\omega)$ with rational dilation factor $M$ is given as formula (4), then the set $\left\{\psi_{j,n}^m(x) = M^{-\frac{j}{2}}\psi^m(M^{-j}x - nq)\right\}_{j,n \in \mathbb{Z},\ 1 \leq m \leq p-q}$ is an orthonormal wavelet basis for $L^2(R)$.

The proof is given as follows.

### 1) The proof of orthogonality

According to the Parseval formula, $\langle \psi_{j,n}^{m_1}(x), \psi_{l,k}^{m_2}(x)\rangle = \frac{1}{2\pi}\langle \hat{\psi}_{j,n}^{m_1}(\omega), \hat{\psi}_{l,k}^{m_2}(\omega)\rangle$, so

$$\langle \psi_{j,n}^{m_1}(x), \psi_{l,k}^{m_2}(x)\rangle$$
$$= \int_{-\infty}^{+\infty} M^{-\frac{j}{2}}\psi^{m_1}(M^{-j}x - nq) \overline{M^{-\frac{l}{2}}\psi^{m_2}(M^{-l}x - kq)}\ dx$$

$$= M^{-\frac{j+l}{2}}\int_{-\infty}^{+\infty} \psi^{m_1}(M^{-j}x - nq) \overline{\psi^{m_2}(M^{-l}x - kq)}\ dx$$

$$= M^{-\frac{j+l}{2}}\frac{1}{2\pi}\int_{-\infty}^{+\infty} M^j\hat{\psi}^{m_1}(M^j\omega)e^{-iM^jnq\omega} \overline{M^l\hat{\psi}^{m_2}(M^l\omega)e^{-iM^lkq\omega}}\ d\omega$$

$$= M^{\frac{j+l}{2}}\frac{1}{2\pi}\int_{-\infty}^{+\infty} \hat{\psi}^{m_1}(M^j\omega) \overline{\hat{\psi}^{m_2}(M^l\omega)}e^{-i(M^jn - M^lk)q\omega}\ d\omega \quad (7)$$

For the formula (7), several cases are discussed in the following.

**A)** $j \neq l$

From formula (5), the support of $\hat{\psi}^{m_1}(M^j\omega)$ is $\left(-M^{-j}(q+m_1)\frac{\pi}{q},\ -M^{-j}(q+m_1-1)\frac{\pi}{q}\right] \cup \left[M^{-j}(q+m_1-1)\frac{\pi}{q},\ M^{-j}(q+m_1)\frac{\pi}{q}\right)$, since $1 \leq m_1 \leq p-q, m_1 \in \mathbb{Z}, p, q \in \mathbb{Z}^+, p > q$, and $M = \frac{p}{q}$, so it is a subset of set $(-M^{-j+1}\pi,\ -M^{-j}\pi] \cup [M^{-j}\pi, M^{-j+1}\pi)$. In the same way, the support of $\hat{\psi}^{m_2}(M^l\omega)$ is $\left(-M^{-l}(q+m_2)\frac{\pi}{q},\ -M^{-l}(q+m_2-1)\frac{\pi}{q}\right] \cup \left[M^{-l}(q+m_2-1)\frac{\pi}{q},\ M^{-l}(q+m_2)\frac{\pi}{q}\right)$, and is a subset of set $(-M^{-l+1}\pi,\ -M^{-l}\pi] \cup [M^{-l}\pi, M^{-l+1}\pi)$. Since $j \neq l$, so the intersection of

$(-M^{-j+1}\pi, -M^{-j}\pi] \cup [M^{-j}\pi, M^{-j+1}\pi)$ and $(-M^{-l+1}\pi, -M^{-l}\pi] \cup [M^{-l}\pi, M^{-l+1}\pi)$ is empty set, then the intersection of the supports of $\hat{\psi}^{m_1}(M^j\omega)$ and $\hat{\psi}^{m_2}(M^l\omega)$ is empty set, so $\langle \psi_{j,n}^{m_1}(x), \psi_{l,k}^{m_2}(x) \rangle = 0$, for $j \neq l$, $j \in \mathbb{Z}$, $l \in \mathbb{Z}$, $1 \leq m_1 \leq p-q$, $1 \leq m_2 \leq p-q$, $m_1 \in \mathbb{Z}$, $m_2 \in \mathbb{Z}$, and any $k, n \in \mathbb{Z}$.

B) $m_1 \neq m_2$

If $m_1 \neq m_2, j \neq l$, then $\langle \psi_{j,n}^{m_1}(x), \psi_{l,k}^{m_2}(x) \rangle = 0$, the proof is given in the above section A.

For $m_1 \neq m_2$, $j = l$, without loss of generality, let $m_2 > m_1$, since $m_1 \in \mathbb{Z}, m_2 \in \mathbb{Z}$, so $m_2 - 1 \geq m_1$, and $p, q \in \mathbb{Z}^+$, $M > 1$, then

$M^{-l}(m_2 + q - 1)\dfrac{\pi}{q} \geq M^{-l}(m_1 + q)\dfrac{\pi}{q}$, so the intersection of the supports of $\hat{\psi}^{m_1}(M^j\omega)$ and $\hat{\psi}^{m_2}(M^l\omega)$ is empty set, then formula (7) is equal to zero.

So, for $m_1 \neq m_2, 1 \leq m_1 \leq p-q$, $1 \leq m_2 \leq p-q$, $m_1, m_2 \in \mathbb{Z}$, and any $k, n \in \mathbb{Z}$, $j, l \in \mathbb{Z}$.

C) $k \neq n$

If $m_1 \neq m_2$, then $\langle \psi_{j,n}^{m_1}(x), \psi_{l,k}^{m_2}(x) \rangle = 0$ for any $j, l \in \mathbb{Z}$, $k, n \in \mathbb{Z}$, the proof is given in the above section A and section B.

If $m_1 = m_2$, $j \neq l$, then $\langle \psi_{j,n}^{m_1}(x), \psi_{l,k}^{m_2}(x) \rangle = 0$ for any $k, n \in \mathbb{Z}$, the proof is given in the above section A.

For $m_1 = m_2$, $j = l$, $k \neq n$, from formula (7),

$$\langle \psi_{j,n}^{m_1}(x), \psi_{l,k}^{m_2}(x) \rangle$$

$$= \langle \psi_{j,n}^{m_1}(x), \psi_{j,k}^{m_1}(x) \rangle$$

$$= M^j \frac{1}{2\pi} \int_{-\infty}^{+\infty} \hat{\psi}^{m_1}(M^j\omega) \overline{\hat{\psi}^{m_1}(M^j\omega)} e^{-i(M^j n - M^j k)q\omega} \, d\omega$$

$$= \frac{M^j}{2\pi} \left[ \int_{M^{-j}(m_1+q-1)\frac{\pi}{q}}^{M^{-j}(m_1+q)\frac{\pi}{q}} (\sqrt{q})^2 e^{-i(n-k)M^j q\omega} \, d\omega \right.$$

$$\left. + \int_{-M^{-j}(m_1+q)\frac{\pi}{q}}^{-M^{-j}(m_1+q-1)\frac{\pi}{q}} (\sqrt{q})^2 e^{-i(n-k)M^j q\omega} \, d\omega \right]$$

$$= q \frac{M^j}{2\pi} 2 \int_{M^{-j}(m_1+q-1)\frac{\pi}{q}}^{M^{-j}(m_1+q)\frac{\pi}{q}} \cos(n-k)M^j q\omega \, d\omega$$

$$= \frac{1}{\pi(n-k)} [\sin(n-k)(q+m_1)\pi - \sin(n-k)(q+m_1-1)\pi] \qquad (8)$$

Since $k \neq n$, $k, n \in \mathbb{Z}$, $q \in \mathbb{Z}^+$, $m_1 \in \mathbb{Z}$, so $(n-k)(q+m_1) \in \mathbb{Z}$, $(n-k)(q+m_1-1) \in \mathbb{Z}$, then

$\frac{1}{\pi(n-k)}[\sin(n-k)(q+m_1)\pi - \sin(n-k)(q+m_1-1)\pi]$ is equal to zero.

So, $\langle \psi_{j,n}^{m_1}(x), \psi_{l,k}^{m_2}(x)\rangle = 0$ for $k \neq n$, $k,n \in \mathbb{Z}$, $1 \leq m_1 \leq p-q$, $1 \leq m_2 \leq p-q$, $m_1, m_2 \in \mathbb{Z}$, and any $j, l \in \mathbb{Z}$.

2) **The proof of normalization**

$$\|\psi_{j,n}^m(x)\|^2$$

$$= \langle \psi_{j,n}^m(x), \psi_{j,n}^m(x)\rangle$$

$$= \frac{1}{2\pi} \langle \hat{\psi}_{j,n}^m(\omega), \hat{\psi}_{j,n}^m(\omega)\rangle$$

$$= \frac{1}{2\pi} \int_{-\infty}^{+\infty} M^{\frac{j}{2}} \hat{\psi}^m(M^j \omega) e^{-inM^j q\omega} \overline{M^{\frac{j}{2}} \hat{\psi}^m(M^j \omega) e^{-inM^j q\omega}} \, d\omega$$

$$= \frac{M^j}{2\pi}\left[\int_{M^{-j}(m+q-1)\frac{\pi}{q}}^{M^{-j}(m+q)\frac{\pi}{q}} (\sqrt{q})^2 \, d\omega + \int_{-M^{-j}(m+q)\frac{\pi}{q}}^{-M^{-j}(m+q-1)\frac{\pi}{q}} (\sqrt{q})^2 \, d\omega\right]$$

$$= \frac{M^j}{2\pi}\left[M^{-j}\frac{\pi}{q}q + M^{-j}\frac{\pi}{q}q\right]$$

$$= 1 \qquad (9)$$

In conclusion,

$\langle \psi_{j,n}^{m_1}(x), \psi_{l,k}^{m_2}(x)\rangle = \delta(j-l)\delta(m_1-m_2)\delta(n-k)$. So, the family $\{\psi_{j,n}^m(x)\}_{j,n \in \mathbb{Z}, 1 \leq m \leq p-q}$ is orthonormal.

3) **The proof of completeness**

Let $W_j^m$ be the collection of functions that have Fourier transforms supported in $\left(-M^{-j}(q+m)\frac{\pi}{q}, -M^{-j}(q+m-1)\frac{\pi}{q}\right] \cup \left[M^{-j}(q+m-1)\frac{\pi}{q}, M^{-j}(q+m)\frac{\pi}{q}\right)$, that is, $M^{-j}(q+m-1)\frac{\pi}{q} \leq |\omega| < M^{-j}(q+m)\frac{\pi}{q}$. Then, if fix $j$, $\bigcup W_j^m$ ($1 \leq m \leq p-q$) is the collection of functions that have Fourier transforms supported in $(-M^{-j+1}\pi, -M^{-j}\pi] \cup [M^{-j}\pi, M^{-j+1}\pi)$. If $j$ varies from $-\infty$ to $+\infty$, then $\bigcup W_j^m = L^2(R)$, $j \in \mathbb{Z}$, $1 \leq m \leq p-q$.

Let $f(x) \in W_j^m$, $\hat{f}(\omega)$ is the Fourier transform of $f(x)$, then the support of $\hat{f}(\omega)$ is $\left(-M^{-j}(q+m)\frac{\pi}{q}, -M^{-j}(q+m-1)\frac{\pi}{q}\right] \cup \left[M^{-j}(q+m-1)\frac{\pi}{q}, M^{-j}(q+m)\frac{\pi}{q}\right)$, so $f(x)$ is a bandpass signal. From the expression of $\hat{\psi}^m(M^j \omega)$ (as shown in (5)), $\frac{T}{\sqrt{q}} \hat{\psi}^m(M^j \omega)$ is a bandpass frequency gating function. Here T is sampling period.

The inverse Fourier transform of $\frac{T}{\sqrt{q}}\hat{\psi}^m(M^j\omega)$ is $\frac{T}{\sqrt{q}}M^{-j}\psi^m(M^{-j}x)$. Based on the reconstruction formula of bandpass signal [12], (10) is got.

$$f(x)=\sum_{n=-\infty}^{\infty}f(nT)\frac{T}{\sqrt{q}}M^{-j}\psi^m[M^{-j}(x-nT)] \qquad (10)$$

$f(x)\in W_j^m$, so, it is a bandpass signal. The higher frequency of the band is $F_H=M^{-j}(q+m)\frac{\pi}{q}$, the lower frequency of the band is $F_L=M^{-j}(q+m-1)\frac{\pi}{q}$, so the bandwidth B is $M^{-j}\frac{\pi}{q}$. Since $\frac{F_H}{B}=q+m$, $q+m$ are integers, so according to the sampling theorem for bandpass signal [12], choosing the sampling frequency $F_S=2B$. Here, $F_S$ is angular frequency, the corresponding natural frequency is $f_S=\frac{F_S}{2\pi}=\frac{2M^{-j}\pi}{q}\frac{1}{2\pi}=\frac{M^{-j}}{q}$, so the sampling period $T=\frac{1}{f_S}=qM^j$. Substitute $T=qM^j$ into formula (10), so

$$f(x)=\sum_{n=-\infty}^{\infty}f(nT)\frac{T}{\sqrt{q}}M^{-j}\psi^m[M^{-j}(x-nT)]$$

$$=\sum_{n=-\infty}^{\infty}f(nT)\frac{qM^j}{\sqrt{q}}M^{-j}\psi^m[M^{-j}(x-nqM^j)]$$

$$=\sum_{n=-\infty}^{\infty}\left(f(nT)M^{\frac{j}{2}}\sqrt{q}\right)M^{-\frac{j}{2}}\psi^m(M^{-j}x-nq) \qquad (11)$$

According to formula (11), $f(nT)M^{\frac{j}{2}}\sqrt{q}$ is a constant, so the set $\left\{M^{-\frac{j}{2}}\psi^m(M^{-j}x-nq)\right\}_{n\in\mathbb{Z}}$ is a basis for $W_j^m$. In the previous section, it has been proved that $M^{-\frac{j}{2}}\psi^m(M^{-j}x-nq)$ is orthonormal, so the set $\left\{M^{-\frac{j}{2}}\psi^m(M^{-j}x-nq)\right\}_{n\in\mathbb{Z}}$ is an orthonormal basis for $W_j^m$. Since $\bigcup W_j^m=L^2(R)$, $j\in\mathbb{Z}$, $1\leq m\leq p-q$, so the set $\left\{\psi_{j,n}^m(x)=M^{-\frac{j}{2}}\psi^m(M^{-j}x-nq)\right\}_{j,n\in\mathbb{Z},\ 1\leq m\leq p-q}$ is an orthonormal wavelet basis for $L^2(R)$. □

**Corollary:** The basis of rational Littlewood-Paley wavelet with rational dilation factor M=p/q first proposed by Auscher is orthonormal if and only if q equals to 1.

**Proof:** Comparing (1) and (3), it can be found that Auscher's rational Littlewood-Paley wavelet is the same as the new Littlewood-Paley wavelet in this letter if the q in (3) is equal to 1. Based on Theorem 1 and 2, the Corollary is proved. □

## V. Conclusion

The orthonormal basis of rational Littlewood-Paley wavelet with rational dilation factor M=p/q first proposed by Auscher does not hold for all rational numbers. In other words, it is not an orthonormal basis if the rational dilation factor M is not an integer. It is orthonormal if and only if q equals to 1. In this letter, the new rational Littlewood-Paley wavelet which holds for all rational numbers is given in Fourier domain and time domain, respectively. In fact, the rational Littlewood-Paley wavelet proposed by Auscher is a special case of the rational Littlewood-Paley wavelet proposed in this letter when q is equal to 1. By rigorous proof, it is shown that the new rational Littlewood-Paley wavelet family $\left\{\psi_{j,n}^{m}(x)=M^{-\frac{j}{2}}\psi^{m}(M^{-j}x-nq)\right\}_{j,n\in\mathbb{Z},\ 1\leq m\leq p-q}$ is an orthonormal wavelet basis for $L^2(R)$.

In future work, based on our previous work [8], the new analysis and synthesis fast algorithm will be studied for rational wavelet transform, and the new rational Littlewood-Paley wavelet will be used to verify the validity of the new analysis and synthesis fast algorithm.